\UseRawInputEncoding
\documentclass[10pt]{article}
\usepackage{relsize}
\usepackage{url}
\usepackage[dvips]{color}
\usepackage{epsfig}
\usepackage{float,amsthm,amssymb,amsfonts}
\usepackage{ amssymb,amsmath,graphicx, amsfonts, latexsym}
\usepackage{xcolor}
\begin{document}
\theoremstyle{plain}
\newtheorem{theorem}{{\bf Theorem}}[section]
\newtheorem{corollary}[theorem]{Corollary}
\newtheorem{lemma}[theorem]{Lemma}
\newtheorem{proposition}[theorem]{Proposition}
\newtheorem{remark}[theorem]{Remark}

\theoremstyle{definition}
\newtheorem{defn}{Definition}
\newtheorem{definition}[theorem]{Definition}
\newtheorem{example}[theorem]{Example}
\newtheorem{conjecture}[theorem]{Conjecture}

\def\im{\mathop{\rm Im}\nolimits}
\def\dom{\mathop{\rm Dom}\nolimits}
\def\rank{\mathop{\rm rank}\nolimits}
\def\nullset{\mbox{\O}}
\def\ker{\mathop{\rm ker}\nolimits}
\def\implies{\; \Longrightarrow \;}

\def\GR{{\cal R}}
\def\GL{{\cal L}}
\def\GH{{\cal H}}
\def\GD{{\cal D}}
\def\GJ{{\cal J}}

\def\set#1{\{ #1\} }
\def\z{\set{0}}
\def\Sing{{\rm Sing}_n}
\def\nullset{\mbox{\O}}

\title{On   injective partial Catalan monoids}
\author{\bf  F. S. Al-Kharousi, A. Umar and   M. M. Zubairu\footnote{Corresponding Author. ~~Email: \emph{mmzubairu.mth@buk.edu.ng}}    \\
\it\small  Department of Mathematics,\\
\it\small College of Science,\\
\it\small Sultan Qaboos University, Oman.\\
\it\small \texttt{fatma9@squ.edu.om} \\[3mm]
\it\small Department of Mathematical Sciences,\\
\it\small Khalifa University, P. O. Box 127788, Sas al Nakhl, Abu Dhabi, UAE\\
\it\small  \texttt{abdullahi.umar@ku.ac.ae}\\[3mm]
\it\small  Department of Mathematical Science,\\
\it\small College of Physical Sciences,\\
\it\small Bayero  University, Kano, Nigeria.\\
\it\small  \texttt{mmzubairu.mth@buk.edu.ng}\\[2mm]}
\date{\today}
\maketitle\

\begin{abstract}
 Let $[n]$ be a finite chain $\{1, 2, \ldots, n\}$, and let $\mathcal{IC}_{n}$ be the semigroup consisting of all isotone and order-decreasing injective partial transformations on $[n]$. In addition, let $\mathcal{Q}^{\prime}_{n} = \{\alpha \in \mathcal{IC}_{n} : \, 1\not \in \textnormal{Dom } \alpha\}$ be the subsemigroup of $\mathcal{IC}_{n}$, consisting of all transformations in $\mathcal{IC}_{n}$, each of whose domains does not contain $1$. For $1 \leq p \leq n$, let
$K(n,p) = \{\alpha \in \mathcal{IC}_{n} : \, |\im \, \alpha| \leq  p\}$
and
$M(n,p) = \{\alpha \in \mathcal{Q}^{\prime}_{n} : \, |\im \, \alpha| \leq p\}$
be the two-sided ideals of $\mathcal{IC}_{n}$ and $\mathcal{Q}^{\prime}_{n}$, respectively. Moreover, let
${RIC}_{p}(n)$ and ${RQ}^{\prime}_{p}(n)$ denote the Rees quotients of $K(n,p)$ and $M(n,p)$, respectively. It is shown in this article that for any \( S \in \{ \mathcal{RIC}_{p}(n), K(n,p) \} \), \( S \) is abundant; \( \mathcal{IC}_{n} \) is ample; and for any \( S \in \{ \mathcal{Q}^{\prime}_{n}, \mathcal{RQ}^{\prime}_{p}(n), M(n,p) \} \), \( S \) is right abundant for all values of \( n \), but not left abundant for \( n \geq 2 \). Furthermore, the ranks of the Rees quotients  ${RIC}_{p}(n)$ and ${RQ}^{\prime}_{p}(n)$  are shown to be equal to the ranks of the two-sided ideals $K(n,p)$ and $M(n,p)$, respectively. These ranks are found to be $\binom{n}{p}+(n-1)\binom{n-2}{p-1}$ and $\binom{n}{p}+(n-2)\binom{n-3}{p-1}$,  respectively. In addition, the ranks of the semigroups  $\mathcal{IC}_{n}$ and $\mathcal{Q}^{\prime}_{n}$ were found to be $2n$ and   $n^{2}-3n+4$, respectively. Finally, we characterize all the maximal subsemigroups of $\mathcal{IC}_{n}$  and $\mathcal{Q}^{\prime}_{n}$. \end{abstract}

\emph{2020 Mathematics Subject Classification. 20M20.}\\
\textbf{Keywords:} Isotone maps, Order decreasing, ample semigroup, Rank properties

\section{Introduction and Preliminaries}

 If $ n $ is a natural number, then $ [n] $ is defined to be the chain $ \{1, 2, \ldots, n\} $. A map $ \alpha $ is called a \emph{partial transformation} if its domain and range are subsets of $ [n] $.  The collection of all partial transformations on $[n]$ under composition of functions is known as \emph{the semigroup of all partial transformations}, which is usually denoted by $\mathcal{P}_{n}$. The collection of all injective partial transformations on $[n]$ is known as  \emph{the semigroup of all partial one-one transformations on} $[n]$ (more commonly
known as \emph{the symmetric inverse semigroup}), and is usually denoted by $\mathcal{I}_{n}$. A transformation $\alpha\in \mathcal{P}_{n}$ is said to be an \emph{ isotone} map (resp., an \emph{anti-tone} map) if  (for all $x,y \in \dom\,\alpha$) $x\leq y$ implies $x\alpha\leq y\alpha$ (resp., $x\alpha\geq y\alpha$); \emph{order decreasing} if (for all $x\in \dom\,\alpha$) $x\alpha\leq x$. The notations  $\mathcal{DI}_n$ and $\mathcal{OI}_n$  shall denote  \emph{the semigroup of all order-decreasing injective partial transformations} on $[n]$ and  \emph{the semigroup of all isotone injective partial transformations} on $[n]$, respectively.

 The notation \( \mathcal{IC}_{n} \) was first introduced in \cite{gmv} to denote the \emph{semigroup of all isotone and order-decreasing injective partial transformations} on \( [n] \). The order of \( \mathcal{IC}_{n} \) was obtained in [\cite{gmv}, Theorem 2.4.8(ii)], and it is known as \( c_{n+1} \), where \( c_{n} = \frac{1}{n} \binom{2n}{n-1} \) is the \emph{$n$-th Catalan number}. Hence, since $\mathcal{IC}_{n}$ is a monoid, we shall refer to it   as the \emph{injective partial Catalan monoid}, which is  obtained from the intersection: \begin{equation}\label{qn111}\mathcal{IC}_{n}= \mathcal{OI}_n\cap \mathcal{DI}_n  .\end{equation} \noindent The monoids $\mathcal{IC}_{n}$, $\mathcal{OI}_n$ and $\mathcal{DI}_n$ have been extensively studied in various contexts, see for example \cite{gmv, g2, lu2015, umar, abu, ua1, ua3, pv}. The composition of two transformations $\alpha $ and $\beta$ in $\mathcal{I}_{n}$ is defined as $x(\alpha\circ\gamma)=((x)\alpha)\gamma$ for all $x\in\dom\, \alpha$. To avoid confusion, we will use the notation $\alpha\beta$ to represent $\alpha\circ\beta$. Moreover, we shall also adopt the notations $F(\alpha)$, $f(\alpha)$, $\text{id}_{[n]}$, $\im \alpha$, $\dom \alpha$, $h(\alpha)$ to mean the set of fixed points of $\alpha$, the number of fixed points of $\alpha$ (i.e., $|F(\alpha)|$),  identity map on $[n]$, the image set of a map $\alpha$, the domain set of a map $\alpha$ and the height of $\alpha$(i.e., $|\im \, \alpha|$), respectively. Moreover, let $s(\alpha)$ be the cardinal of the set $S(\alpha)=\{x\in[n]: \, x\alpha\neq x\}$, $s(\alpha)$ is called the \emph{shift} of $\alpha$.
 Now let
\begin{equation}\label{qn1}
 	\mathcal{Q}_{n} = \{\alpha \in \mathcal{IC}_{n} : 1 \in \textnormal{Dom } \alpha \}
 \end{equation}
\noindent be the set of all maps in $\mathcal{IC}_{n}$ each of whose domains contain  1 and
\begin{equation}\label{qn2}
	\mathcal{Q}^{\prime}_n = \{\alpha \in \mathcal{IC}_{n} : 1 \notin \text{Dom } \alpha\}
\end{equation}
 \noindent be the set of all maps in $\mathcal{IC}_{n}$ each of whose domains do not contain  1. In other words, $\mathcal{Q}^{\prime}_n$ is the set complement of $\mathcal{Q}_{n}$. Notice that for $\alpha$, $\beta\in \mathcal{Q}_{n}$, since $1\in \dom \, \alpha$ then $1\alpha\beta=(1\alpha)\beta=1\beta=1$, which shows that $1\in \dom \, \alpha\beta$; and for any $x\in \dom \, \alpha$, $x\alpha\beta\leq x\beta\leq x$. This ensures that $\mathcal{Q}_{n}$ is a subsemigroup of $\mathcal{IC}_{n}$. Similarly, one can demonstrate easily that $\mathcal{Q}^{\prime}_n$ is a subsemigroup of $\mathcal{IC}_{n}$. It is worth noting that the semigroup $\mathcal{Q}_{n}$ contains the identity element, while $\mathcal{Q}^{\prime}_n$ does not. Furthermore, since for any $\alpha\in \mathcal{Q}_n$, $1\in \dom \, \alpha$, it follows that 1 is a fixed point, and  so, it can easily be demonstrated that $\mathcal{Q}_n$ (on $n$ elements) is isomorphic to the injective partial Catalan monoid $\mathcal{IC}_{n-1}$ (on $n-1$ elements), and so, the order of $\mathcal{Q}_n$ is: \[c_{n}=\frac{1}{n}\binom{2n}{n-1}.\] \noindent This number corresponds to the $n$'th Catalan number, therefore  $\mathcal{Q}_n$ is the \emph{injective  partial Catalan monoid} on $[n-1]$ chain.  Therefore, the algebraic, combinatorial and rank properties  of $\mathcal{Q}_n$ can easily be obtain from that of $\mathcal{IC}_{n-1}$ by isomorphism. Hence we need not study these properties on the semigroup $\mathcal{Q}_n$. Now one can easily see that the order of $\mathcal{Q}^{\prime}_n$ is \[t_{n}=c_{n+1}-c_{n}=\frac{3}{n+2}\binom{2n}{n-1}=\frac{3nc_{n}}{n+2}.\]
 \noindent The number $t_{n}$ corresponds to the sequence \underline{$A000245$}  in The On-Line Encyclopedia of Integer Sequences \cite{slone}. There are certain algebraic and combinatorial properties of the semigroup $\mathcal{IC}_{n}$ that have been explored; see, for example, \cite{auc, abu}.

 Now, for $1\le p\le n$,  let \begin{equation} \label{kn} K(n,p)=\{\alpha\in  \mathcal{IC}_{n}: \, |\im \, \alpha|\le p\}\end{equation}
   \noindent and
   \begin{equation}\label{mn} M(n,p)=\{\alpha\in  \mathcal{Q}^{\prime}_{n}: \, |\im \, \alpha|\le p\}\end{equation}

  \noindent be the two-sided ideals of $\mathcal{IC}_{n}$ and $\mathcal{Q}^{\prime}_{n}$, respectively, consisting of all decreasing isotone injective partial maps with a height of no more than $p$.

  Furthermore, for $p\geq 1$,  let \begin{equation}\label{knn} {RIC}_{n}(p)= K(n,p)/ K(n, p-1);  \end{equation}
\noindent and for $p\geq 2$

     \begin{equation}\label{mnnn} {RQ}^{\prime}_{n}(p)= M(n,p)/M(n, p-1)  \end{equation}

  \noindent be the Rees quotient semigroups of $K(n,p)$ and $M(n,p)$, respectively. The elements of ${RIC}_{n}(p)$ (resp., ${RQ}^{\prime}_{n}(p)$) can be considered as the elements of $\mathcal{IC}_{n}$ (resp., $\mathcal{Q}^{\prime}_{n}$) with exactly height $p$. The product of two elements of ${RIC}_{n}(p)$ (resp.,  ${RQ}^{\prime}_{n}$) is 0 if their product in ${RIC}_{n}(p)$ (resp., ${RQ}^{\prime}_{n}(p)$) has a height strictly less than $p$; otherwise it is in ${RIC}_{n}(p)$ (resp.,${RQ}^{\prime}_{n}(p)$).

   Algebraic and rank  properties of these semigroups have not been explored to our knowledge. In this paper, we are going to study certain algebraic and rank properties of these  semigroups. In Section 2, we characterize Green's relations and their starred analogues for all semigroup \( S \in \{ \mathcal{IC}_{n}, \mathcal{Q}^{\prime}_{n}, K(n,p), M(n,p), RIC_{n}(p), RQ_{n}(p) \} \) and show that for any \( S \in \{ \mathcal{RIC}_{p}(n), K(n,p) \} \), \( S \) is abundant. Furthermore, we show that \( \mathcal{IC}_{n} \) is ample, and for any \( S \in \{ \mathcal{Q}^{\prime}_{n}, \mathcal{RQ}^{\prime}_{p}(n), M(n,p) \} \), \( S \) is right abundant for all values of \( n \), but not left abundant for \( n \geq 2\). Moreover, the rank of every \( S \in \{ \mathcal{IC}_{n}, \mathcal{Q}^{\prime}_{n},  M(n,p), RIC_{n}(p), RQ_{n}(p) \}  \) is completely determined in Section 3 through the concept of quasi-idempotents introduced in \cite{al2} and the concept of essential elements introduced in that section. It is worth noting that the rank of $K(n,r)$ has been investigated in \cite{lu}. Finally, we dedicate Section 4 to characterizing all the maximal subsemigroups of the monoid \( \mathcal{IC}_{n} \) and the semigroup \( \mathcal{Q}^{\prime}_{n} \).
 For more details about   terms and concepts in semigroup theory, we recommend that the reader consult  the books by  Howie \cite{howi} and Higgins \cite{ph}.

\indent We shall represent every $\alpha\in \mathcal{IC}_{n}$ in two line notation as:
\begin{equation}\label{1}\alpha=\begin{pmatrix}x_1&\ldots&x_p\\a_1&\ldots&a_p\end{pmatrix}   \quad  (1\le p\le n),\end{equation}  where  $a_{i}\leq x_{i}$ for all $1\leq i\leq p$. Moreover, we may, without loss of generality, assume that $1\leq x_{1}<\ldots<x_{p}\leq n$ and  $1\leq a_{1}<a_{2}<\ldots<a_{p}\leq n$, since $\alpha$ is an isotone and injective map.

\section{Green's relations and their starred analogue}

In a semigroup \( S \), an element \( a \in S \) is called \emph{regular} if there exists an element \( b \in S \) such that \( a = aba \). Additionally, a semigroup \( S \) is referred to as a \emph{regular semigroup} if every element in \( S \) is regular. In semigroup theory, there are five Green's relations, specifically denoted as \( \mathcal{L}, \mathcal{R}, \mathcal{D}, \mathcal{J}, \) and \( \mathcal{H} \). These relations are defined in terms of the principal (left, right, or two-sided) ideals generated by the elements in \( S \). For a detailed explanation of these relations, we suggest that the reader consult Howie \cite{howi}. As noted in \cite{FOUN}, in case of ambiguity, we denote a relation \( \mathcal{K} \) on \( S \) by \( \mathcal{K}(S) \). It is a well-established fact in finite semigroups that the relations \( \mathcal{D} \) and \( \mathcal{J} \) are equivalent (see \cite{howi}, Proposition 2.1.4). As a result, we will focus on describing the relations \( \mathcal{L}, \mathcal{R},  \mathcal{D}, \) and \( \mathcal{H} \) in the injective partial Catalan monoid \( \mathcal{IC}_{n} \) and its subsemigroup \( \mathcal{Q}^{\prime}_{n}\).

It is well known that the semigroup \( \mathcal{DI}_n \) is \( \mathcal{J}\)-trivial (see [\cite{auc}, Lemma 2.2]), and since \( \mathcal{IC}_{n} \) is a subsemigroup of \( \mathcal{DI}_n \), it follows that for any relation say \( \mathcal{K} \in \{\mathcal{L}, \mathcal{R}, \mathcal{H}, \mathcal{J}\}\), we have
\(
\mathcal{K}(\mathcal{IC}_{n}) =  \mathcal{K}(\mathcal{DI}_n) \cap (\mathcal{IC}_{n} \times \mathcal{IC}_{n}).
\)
Consequently, we have the following result.

\begin{theorem} The semigroup $\mathcal{IC}_{n}$ (resp., $\mathcal{Q}^{\prime}_{n}$) is $\mathcal{J}-$trivial.
 \end{theorem}

  As a consequence of the above theorem, we readily have the following corollary.

 \begin{corollary}\label{rem1} Let $ S\in \{\mathcal{IC}_{n}, \mathcal{Q}^{\prime}_{n}\}$ and let $\alpha \in S$. Then $\alpha$ is regular if and only if $\alpha$ is an idempotent. Hence, the semigroup $S \in \{\mathcal{IC}_{n}, \mathcal{Q}^{\prime}_{n}\}$ is non-regular.
 \end{corollary}
\begin{proof} The result follows from the fact that in an  $\mathcal{R}$-trivial semigroup, every non-idempotent element is not regular.
 \end{proof}

As a consequence of the  above theorem, we deduce the following characterizations of Green's equivalences on the semigroup $S$ in $\{{RIC}_{n}(p), \, {RQ}^{\prime}_{n}(p), \, M(n,p), \, K(n,p) \}$.

 \begin{theorem} Let $S\in \{{RIC}_{n}(p), \, {RQ}^{\prime}_{n}(p), \, M(n,p), \, K(n,p) \}$. Then $S$ is $\mathcal{J}-$trivial.  Hence, for $p \geq 2$, the semigroup $S$ is non-regular.
\end{theorem}

When a semigroup is not regular, it is common practice to examine the starred Green's relations in order to classify its algebraic structural properties. There are five starred Green's equivalences: \( \mathcal{L}^*, \mathcal{R}^*, \mathcal{D}^*, \mathcal{J}^*, \) and \( \mathcal{H}^* \). The relation \( \mathcal{D}^* \) is defined as the join of \( \mathcal{L}^* \) and \( \mathcal{R}^* \), whereas \( \mathcal{H}^* \) represents the intersection of \( \mathcal{L}^* \) and \( \mathcal{R}^* \). A semigroup \( S \) is termed \emph{left abundant} if every \( \mathcal{L}^* \)-class contains an idempotent; it is called \emph{right abundant} if every \( \mathcal{R}^* \)-class contains an idempotent; and it is considered \emph{abundant} if each \( \mathcal{L}^* \)-class and each \( \mathcal{R}^* \)-class of \( S \) includes an idempotent. These classifications of semigroups were introduced by Fountain \cite{FOUN, FOUN2}. For definitions of these relations and their properties, we suggest that the reader consult Fountain \cite{FOUN2}. Therefore, we will concentrate on describing the starred analogues of Green's equivalences on \( S \in \{\mathcal{IC}_{n},  \mathcal{Q}^{\prime}_{n}, {RIC}_{n}(p), {RQ}^{\prime}_{n}(p), M(n,p), K(n,p) \} \).

 Numerous classes of transformation semigroups have been demonstrated to be either left abundant, right abundant, or abundant; see, for instance, \cite{elq, al1, um, umar, quasi, ua3, zuf, zuf2, zm1}. Before we characterize the starred Green's relations, we require the following definitions and lemmas from \cite{elq, quasi}: A subsemigroup $U$ of a semigroup $S$ is called \emph{full} if $E(U)=E(S)$; it is referred to as an \emph{inverse ideal} of $S$ if for every $u \in U$, there exists $u^{\prime} \in S$ such that $uu^{\prime}u = u$ and both $u^{\prime}u$ and $uu^{\prime}$ are in $U$.

 \begin{lemma}[\cite{quasi}, Lemma 3.1.8.]\label{inv1}  Every inverse ideal $U$ of a semigroup $S$ is abundant.
 \end{lemma}

 \begin{lemma} [\cite{quasi}, Lemma 3.1.9.]  \label{inv2}  Let $U$ be an inverse ideal of a semigroup $S$. Then (1)  $\mathcal{L}^{*} (U) = \mathcal{L}(S) \cap (U \times U)$; (2) $\mathcal{R}^{*}( U) = \mathcal{R}(S) \cap(U \times U)$; (3) $\mathcal{H}^{*}( U) = \mathcal{H}(S) \cap (U \times U).$
 \end{lemma}

 We now have the following result.
 \begin{theorem}\label{inv} Let \(\mathcal{IC}_{n}\)  be as defined in \eqref{qn111}. Then \(\mathcal{IC}_{n}\) is an inverse ideal of  $\mathcal{I}_{n}$.
 \end{theorem}
 \begin{proof} Let $\alpha\in \mathcal{IC}_{n}$ be expressed as in \eqref{1}. Then $\alpha^{-1}\in \mathcal{I}_{n}$. Notice that for $1\leq i\leq p$ \[x_{i}\alpha\alpha^{-1}\alpha=a_{i}\alpha^{-1}\alpha=x_{i}\alpha.\]\noindent Moreover, $\alpha\alpha^{-1}=\text{Id}_{\dom \, \alpha}\in \mathcal{IC}_{n}$, and  $\alpha^{-1}\alpha=\text{Id}_{\im \, \alpha}\in \mathcal{IC}_{n}$. This shows that \(\mathcal{IC}_{n}\) is an inverse ideal of \(\mathcal{I}_{n}\), as required.
 \end{proof}
 \begin{remark}\label{gg} \begin{itemize} \item[(i)]  It is easy to see that $\mathcal{IC}_{n}$ is a full subsemigroup of $\mathcal{I}_{n}^{-}$ (where $\mathcal{I}_{n}^{-}$ denotes the semigroup of all decreasing partial injective transformations on $[n]$, see \cite{ua2}) in the sense that $E(\mathcal{IC}_{n}) = E(\mathcal{I}_{n})$. Hence, $E(\mathcal{IC}_{n})$ is a \emph{semilattice} (by semilattice, we mean a commutative semigroup such that all its elements are idempotents).

 \item[(ii)] The semigroup $\mathcal{Q}^{\prime}_{n}$ is a right inverse ideal of $\mathcal{I}_{n}$ but not necessarily a left inverse ideal since $\alpha^{-1}\alpha=\text{id}_{\im \, \alpha}$ is not in $\mathcal{Q}^{\prime}_{n}$ whenever $a_{1}=1$.\end{itemize}
 \end{remark}

 Hence, we have the following results.

\begin{theorem}\label{ab}
	Let $\mathcal{IC}_{n} \ \text{and } \mathcal{Q}^{\prime}_{n}$ be as defined in \eqref{qn111} and \eqref{qn2}, respectively. Then\begin{itemize}
                                                                                           \item[(i)] $\mathcal{IC}_{n}$ is   abundant;
                                                                                                                                      \item[(ii)] $\mathcal{Q}^{\prime}_{n}$ is a right abundant.
                                                                                                                                    \end{itemize}
                                                                                                                                    \end{theorem}

\begin{proof}
	(i)  follows from Theorem \ref{inv}  and Lemma \ref{inv1}, while
(ii)  follows from  Remark \ref{gg} and Lemma \ref{inv1}.
\end{proof}

Now we have the following result.
 \begin{theorem} \label{a1}
 Let $\mathcal{IC}_{n}$ be as defined in \eqref{qn111}. Then  for $\alpha, \beta\in \mathcal{IC}_{n}$, we have:
 \begin{itemize}
   \item[(i)] $\alpha\mathcal{L}^*\beta$  if and only $\im  \alpha = \im  \beta$;
   \item[(ii)] $\alpha\mathcal{R}^*\beta$ if and only if $\dom \,   \alpha = \dom \,  \beta$;
   \item[(iii)] $\alpha\mathcal{H}^*\beta$ if and only if $\alpha=\beta$.
 \end{itemize}
 \end{theorem}

\begin{proof}
\begin{itemize} \item[(i)] and (ii) follow from Theorem \ref{inv}, Lemma \ref{inv2} and the proof of [\cite{ral}, Lemma 2.], while (iii) follows from (i) and (ii) and the fact that $\alpha$ and $\beta$ are isotone.
\end{itemize}
\end{proof}
\begin{remark}\label{a5} The results of the above theorem also hold on $\mathcal{Q}^{\prime}_{n}$. The proof of (i) can be demonstrated using the techniques of that of [\cite{ua1}, Lemma 2.5 (i)]. While (ii) and (iii) are as in that of the above theorem, respectively.
\end{remark}

An abundant semigroup \( S \) is termed \emph{adequate} if \( E(S) \) forms a semilattice. Clearly, inverse semigroups qualify as adequate since, in this scenario, \( \mathcal{L}^* = \mathcal{L} \) and \( \mathcal{R}^* = \mathcal{R} \) \cite{FOUN2}. Thus, by Remark \ref{gg}(i), we see that \( \mathcal{IC}_{n} \) is an adequate semigroup. For an element \( a \) of an adequate semigroup \( S \), the (unique) idempotent in the \( \mathcal{L}^* \)-class (or \( \mathcal{R}^* \)-class) containing \( a \) will be denoted by \( a^{*} \) (or \( a^{+} \)). An adequate semigroup \( S \) is said to be \emph{ample} if \( ea = a(ea)^{*} \) and \( ae = (ae)^{+}a \) for all elements \( a \) in \( S \) and all idempotents \( e \) in \( S \). Now, since \( \mathcal{IC}_{n}\) is a subsemigroup of \( \mathcal{I}_{n}^{-}\), and since \( E(\mathcal{IC}_{n}) = E(\mathcal{I}_{n}^{-})\), and given that \( \mathcal{I}_{n}^{-}\) is an ample semigroup, then by Remark \ref{gg}(i), we readily have the following result.

\begin{theorem} Let $\mathcal{IC}_{n}$ be as defined in \eqref{qn111}. Then $\mathcal{IC}_{n}$ is an ample semigroup for all $n\geq 2$.
\end{theorem}
Now observe that $E(\mathcal{Q}^{\prime}_{n})=\{\epsilon\in E(\mathcal{IC}_{n}): \, 1\notin \dom \epsilon \}$, and so, it is clear that for any $\epsilon_{1}, \epsilon_{2}\in  E(\mathcal{Q}^{\prime}_{n})$, $1\notin \dom \epsilon_{1}\epsilon_{2}$, and as such $\epsilon_{1}\epsilon_{2}\in E(\mathcal{Q}^{\prime}_{n})$. Thus, $E(\mathcal{Q}^{\prime}_{n})$ is a semilattice. Therefore by Theorem \ref{ab}(ii), we see that  $\mathcal{Q}^{\prime}_{n}$ is  right adequate. Moreover, since $E(\mathcal{Q}^{\prime}_{n})\subset E(\mathcal{IC}_{n})$, then we readily have the following result.
\begin{theorem} Let $\mathcal{Q}^{\prime}_{n}$ be as defined in \eqref{qn2}. Then $\mathcal{Q}^{\prime}_{n}$ is a right ample semigroup for all $n$.
\end{theorem}

It is important to note that on the semigroups $\mathcal{IC}_{n} \ \text{and } \mathcal{Q}^{\prime}_{n}$, the relations $\mathcal{L}^{*}$ and $\mathcal{R}^{*}$ do not commute as in the lemma below.
\begin{lemma}\label{commute}

 \begin{itemize}\item[(i)] For $n\geq 2$, on the injective partial Catalan monoid $\mathcal{IC}_{n}$, we have $\mathcal{R}^{*}\circ\mathcal{L}^{*} \neq \mathcal{L}^{*}\circ\mathcal{R}^{*}.$
\item[(ii)]  For $n\geq 3$, on the semigroup $ \mathcal{Q}^{\prime}_{n}$, we have $\mathcal{R}^{*}\circ\mathcal{L}^{*} \neq \mathcal{L}^{*}\circ\mathcal{R}^{*}.$\end{itemize}
\end{lemma}
\begin{proof} In (i), take \[\alpha=\left(\begin{array}{c}
                                                                            1   \\
                                                                            {1}
                                                                          \end{array}
\right) \text{ and } \beta=\left(\begin{array}{c}
                                                                            2   \\
                                                                            {2}
                                                                          \end{array}
\right) \] \noindent in $\mathcal{IC}_{n}$; and in (ii) take \[\alpha=\left(\begin{array}{c}
                                                                              2 \\
                                                                            2
                                                                          \end{array}
\right) \text{ and } \beta=\left(\begin{array}{cc}
                                                                             3 \\
                                                                            3
                                                                          \end{array}
\right)\in \mathcal{Q}^{\prime}_{n}. \]
\noindent Clearly, in either case (by Theorem \ref{a1}), we see that \( (\alpha, \, \beta) \in \mathcal{L}^{*} \circ \mathcal{R}^{*} \) but \( (\alpha, \beta) \notin \mathcal{R}^{*} \circ \mathcal{L}^{*} \), as required.
\end{proof}

Prior to characterizing the relations $\mathcal{D}^{*}$ and $\mathcal{J}^{*}$, we need the following analogues of [\cite{ua2}, Lemmas 2.7, 2.8, \& 2.9] (respectively), which will be stated without proof.
\begin{lemma}\label{ol} Let $S\in \{\mathcal{IC}_{n}, \mathcal{Q}^{\prime}_{n}\}$ and let $\alpha \in S$ be expressed as in \eqref{1}. Then there exists $\beta\in \mathcal{IC}_{n}$ with $\im \, \beta=\{1,\ldots, p\}$ such that $(\alpha, \, \beta)\in \mathcal{R}^{*}$.
\end{lemma}

\begin{lemma}\label{l} Let $S\in \{\mathcal{IC}_{n}, \mathcal{Q}^{\prime}_{n}\}$, and let $\alpha \in S$. Then there exists $\beta\in \mathcal{IC}_{n}$ with $\dom \, \beta=\{n-p+1,\ldots, n\}$ such that $(\alpha, \, \beta)\in \mathcal{L}^{*}$.
\end{lemma}
\begin{lemma}\label{j} Let $S\in \{\mathcal{IC}_{n}, \mathcal{Q}^{\prime}_{n}\}$, and let $\alpha, \, \beta\in S$. If $\alpha\in J^{*}_{\beta}$, then $|\im \, \alpha|\leq |\im \, \beta|$.
\end{lemma}

On the semigroup \( S \in \{ \mathcal{IC}_{n}, \mathcal{Q}^{\prime}_{n} \} \), we define a relation \( \mathcal{K} \) on \( S \) as follows: \( (\alpha,\,\beta) \in \mathcal{K} \) if and only if \( |\im\,\alpha| = |\im\,\beta| \). It follows that \( \mathcal{D} \subseteq \mathcal{K} \). By applying Theorem \ref{a1} and Lemmas \ref{ol}-\ref{j}, along with the techniques used in the proofs of [\cite{ua1}, Lemmas 2.10 \& 2.14], we can readily derive the following results.

\begin{lemma}\label{uaaaa} On the semigroup $S\in \{\mathcal{IC}_{n}, \mathcal{Q}^{\prime}_{n}\}$, we have $\mathcal{D}^{*}=\mathcal{R}^{*}\circ\mathcal{L}^{*}\circ\mathcal{R}^{*}=\mathcal{L}^{*}\circ\mathcal{R}^{*}\circ\mathcal{L}^{*}$.
\end{lemma}
\begin{corollary} Let $S\in \{\mathcal{IC}_{n}, \mathcal{Q}^{\prime}_{n}\}$ and let $\alpha, \beta \in S$. Then $(\alpha, \, \beta)\in \mathcal{D}^{*}$ if and only if $|\im \, \alpha| = |\im \, \beta|$.
\end{corollary}
\begin{lemma}\label{uaaaaa} On the semigroup $S\in \{\mathcal{IC}_{n}, \mathcal{Q}^{\prime}_{n}\}$, we have $\mathcal{J}^{*}=\mathcal{D}^{*}$.
\end{lemma}

Now we deduce that:
\begin{lemma}\label{uaaa} On the semigroups  ${RIC}_{n}(p)$ and  ${RQ}^{\prime}_{n}$, we have $\mathcal{D}^{*}=\mathcal{R}^{*}\circ\mathcal{L}^{*}\circ\mathcal{R}^{*}=\mathcal{L}^{*}\circ\mathcal{R}^{*}\circ\mathcal{L}^{*}$.
\end{lemma}

\begin{lemma}\label{un} On the semigroup $S$ in $\{\mathcal{IC}_{n}, \, \mathcal{Q}^{\prime}_{n}, \, {RQ}^{\prime}_{n}(p), \, {RIC}_{n}(p),  \, M(n,p), \,  K(n,p) \}$, every $\mathcal{R}^{*}-$class  contains a unique idempotent.
\end{lemma}
\begin{proof} This follows from the fact that  $\dom \, \alpha$ can only admit one image subset of $[n]$ so that $\alpha$ is an idempotent, by the decreasing property of $\alpha$.\end{proof}

\begin{remark}\label{ls}\begin{itemize}
             \item[(i)]  It is now clear from Theorem \ref{a1}(i) and (ii) that, for each \( 1 \leq p \leq n \), the number of \( \mathcal{L}^{*} \)-classes is equal to the number of \( \mathcal{R}^{*} \)-classes in \( J^{*}_{p} = \{ \alpha \in \mathcal{IC}_{n} : |\im\,\alpha| = p \} \), which is equal to the number of all possible subsets of \( [n] \), each of order \( p \); that is, \( \binom{n}{p} \);
                 \item[(ii)] It is also clear from Remark \ref{a5} that, for each \( 0 \leq p \leq n-1 \), the number of \( \mathcal{L}^{*} \)-classes in \( J^{*}_{p} = \{ \alpha \in \mathcal{Q}^{\prime}_{n} : |\im\,\alpha| = p \} \) is \( \binom{n}{p} \), while the number of \( \mathcal{R}^{*} \)-classes is \( \binom{n-1}{p} \) since \( 1 \notin \dom\,\alpha \);
             \item[(iii)] If \( S \in \{ {RIC}_{n}(p), {RQ}^{\prime}_{n}(p), M(n,p), K(n,p) \} \), then the characterizations of the starred Green's relations in Theorem \ref{a1} also hold on \( S \);
                 \item[(iv)] Notice that each idempotent of height $p$ is a partial identity, as such  for each $1\leq p\leq n$, $|E(J^*_{p})|=|E(RIC_{n}(p))\setminus\{0\}|=\binom{n}{p}$.
           \end{itemize}
\end{remark}

Thus, the semigroups \( K(n,p) \) and \( M(n,p) \), like \( \mathcal{IC}_{n} \) and \( \mathcal{Q}^{\prime}_{n} \), respectively, are each the union of the \( \mathcal{J}^{*} \) classes
\(
J_{0}^{*}, J_{1}^{*}, \ldots, J_{p}^{*}
\)
and
\(
J_{1}^{*}, \ldots, J_{p}^{*},
\)
respectively, where
\[
J_{p}^{*} = \{ \alpha \in K(n,p)  ~(\text{resp., } \alpha \in M(n,p)) : |\im\,\alpha| = p \}.
\]

Furthermore, $K(n,p)$ and $M(n, p)$ have  ${\binom{n}{p}}$ and ${\binom{n-1}{p}}$ $\mathcal{R}^{*}-$classes, respectively,  and both $K(n,p)$ and $M(n,p)$ have $\binom{n}{p}$ $\mathcal{L}^{*}-$classes in each $J^{*}_{p}$. Therefore, the Rees quotient semigroups ${RIC}_{n}(p)$ and ${RQ}^{\prime}_{n}(p)$ have ${\binom{n}{p}+1}$ and ${\binom{n-1}{p}+1}$ $\mathcal{R}^{*}-$classes, respectively,  and both have $\binom{n}{p}+1$ $\mathcal{L}^{*}-$classes. (The term 1 is derived from the singleton class containing the zero element in every instance.)

  We now present the following lemmas, which can be easily seen from Remark \ref{ls}.
 \begin{lemma}\label{e2} For all $0\leq p\leq n$, if  $J^{*}_{p}=\{\alpha\in \mathcal{IC}_{n}: \, |\im \, \alpha|=p\}$, then $|E(J^{*}_{p})|={\binom{n}{p}}$.
  \end{lemma}
  \begin{lemma}\label{e3} For all $0\leq p\leq n-1$, if  $J^{*}_{p}=\{\alpha\in \mathcal{Q}^{\prime}_{n}: \, |\im \, \alpha|=p\}$, then $|E(J^{*}_{p})|={\binom{n-1}{p}}$.
  \end{lemma}

  It is a known fact from [\cite{gmv}, Proposition 14.3.1] that $|E(\mathcal{IC}_{n})|=2^{n}$. This result can also be obtained by summing all the idempotents of height within the range $0\leq p\leq n$ in Lemma \ref{e2}. We now present the following result.
  \begin{theorem} Let $\mathcal{Q}^{\prime}_{n}$ be as defined in \eqref{qn1}. Then  $|E(\mathcal{Q}^{\prime}_{n})|=2^{n-1}$.
  \end{theorem}
  \begin{proof}The result follows from Lemma \ref{e3} by summing up the number of idempotents elements of height $p$ (i.e., $\binom{n-1}{p}$) from $p = 0$ to $p = n-1$.
  \end{proof}

\section{Rank properties}
Let $S$ represent a semigroup, and let $A$ be a nonempty subset of $S$. The \emph{smallest subsemigroup} of $S$ that includes $A$ is referred to as the \emph{subsemigroup generated by $A$}, typically denoted by $\langle A \rangle$. If there exists a finite subset $A$ of  $S$ such that $\langle A \rangle$ is equal to $S$, then $S$ is termed a \emph{finitely-generated semigroup}. The \emph{rank} of a finitely generated semigroup $S$ is defined as the smallest size of a subset $A$ for which $\langle A \rangle$ is equal to $S$. Specifically,

$$
\text{rank}(S) = \min\{|A| : \langle A \rangle = S\}.
$$

For a more comprehensive exploration of ranks in semigroup theory, we suggest that the reader consult \cite{hrb, hrb2}. Various authors have investigated the ranks, idempotent ranks, and nilpotent ranks of different classes of semigroups of transformation on finite chain $[n]$. In particular, we highlight the contributions of Garba \cite{g1, g2, gu1}, Gomes and Howie \cite{gm, gm2, gm3}, Howie and McFadden \cite{hf}, Umar \cite{umar, ua1, ua}, and Zubairu \emph{et al.} \cite{zm1, zuf, zuf2}. It is well known that the semigroup $\mathcal{IC}_{n}$ is generated by its idempotent and elements of height $n-1$, see [\cite{gmv}, Theorem 14.4.1] and the rank of it two sided ideals have been investigated in \cite{lu}. However, to the best of our knowledge, the minimum generating set of  the Rees Quotients ${RIC}_{n}(p)$ and ${RQ}^{\prime}_{n}(p)$ of the two-sided ideals $K(n,p)$ and $M(n,p)$ of $\mathcal{IC}_{n}$  and $\mathcal{Q}^{\prime}_{n}$, respectively have not been addressed. The aim of this paper is to compute the rank of the Rees quotients ${RIC}_{n}(p)$ and ${RQ}^{\prime}_{n}(p)$; and the rank of the two-sided ideal  $M(n,p)$ of the  subsemigroup $\mathcal{Q}^{\prime}_{n}$, respectively, thereby determining the rank of $\mathcal{Q}^{\prime}_{n}$ as specific instance. Since $\mathcal{IC}_{n}$ and its subsemigroup $\mathcal{Q}^{\prime}_{n}$ are $\mathcal{J}$-trivial semigroups, then in line with \cite{by}, they admit a minimum generating set.

 We will begin our findings by first introducing the following definition (see [\cite{al2}, Introduction]).
 \begin{definition} An element $\alpha \in \mathcal{IC}_{n}$  is called \emph{quasi-idempotent} if $\alpha^{2}$ is an idempotent. Equivalently,  $\alpha$  in $\mathcal{IC}_{n}$  is \emph{quasi-idempotent} if $\alpha^{4}=\alpha^{2}$.  \end{definition}
 We now note the following.
 \begin{remark}\label{rr} \begin{itemize}\item[(i)]On the semigroup $\mathcal{IC}_{n}$,  every quasi-idempotent of height $p$  and of shift 1 has the form: \begin{equation}\label{r}\varepsilon=\begin{pmatrix}y_1&\cdots&y_{i-1}&y_{i}&y_{i+1}&\cdots&y_p\\y_1&\cdots&y_{i-1}&y_{i}\varepsilon&y_{i+1}&\cdots&y_p\end{pmatrix},\end{equation} \noindent where $1\leq y_{1}<\cdots<y_{i-1}<y_{i}\varepsilon< y_{i}<\cdots<y_{p}\leq n$. Notice also that quasi-idempotents of shift 1 are not idempotents.

 \item[(ii)] Every idempotent element is quasi-idempotent but the converse is not necessarily true.
 \end{itemize}
 \end{remark}

 The following lemma is equivalent to [\cite{gmv}, Lemma 14.4.2]. We also provide an equivalent proof of [\cite{gmv}, Lemma 14.4.2] using the concept of quasi-idempotents.
  \begin{lemma}\label{hq1} The semigroup $\mathcal{IC}_{n}$ is quasi-idempotent  generated.
 \end{lemma}
 \begin{proof} Let $\alpha\in \mathcal{IC}_{n}$ be as expressed in \eqref{1}. Now for $1\le i\leq p$ define \[\varepsilon_{i}=\begin{pmatrix}a_1&\cdots&a_{i-1}&x_{i}&x_{i+1}&\cdots&x_p\\a_1&\cdots&a_{i-1}&a_{i}&x_{i+1}&\cdots&x_p\end{pmatrix} .\] \noindent Clearly, for each $i$, $\varepsilon_{i}\in J^{*}_{p}$. Notice that if $a_{i}=x_{i}$, then $\varepsilon_{i}$ is an idempotent, otherwise it is a quasi-idempotent of shift 1.
 Now,  observe that:
   \begin{align*}\varepsilon_{1}\cdots \varepsilon_{p}=&
\begin{pmatrix}
x_{1} & x_2 & \cdots &  x_p \\
a_1 & x_2 & \cdots &  x_p
\end{pmatrix}
\begin{pmatrix}
a_1 & x_2 & x_{3} & \cdots &  x_{p} \\
a_1 & a_2 & x_{3} & \cdots &  x_{p}
\end{pmatrix}
\cdots
\begin{pmatrix}
a_1 & \cdots & a_{p-2} & x_{p-1} &x_{p}\\
a_1 &  \cdots & a_{p-2} & a_{p-1}& x_{p}
\end{pmatrix}
 \begin{pmatrix}
a_1 & \cdots &  a_{p-1} &x_{p} \\
a_1 &  \cdots &  a_{p-1} & a_{p}
\end{pmatrix}
\\=& \begin{pmatrix}
x_1 &  \cdots &  x_{p} \\
a_1 &  \cdots &  a_{p}
\end{pmatrix}=\alpha,
\end{align*}
\noindent as postulated.
 \end{proof}
\begin{remark} Notice that in the above proof we have shown that $\mathcal{IC}_{n}$ is generated by idempotents and quasi-idempotents of shift 1.
\end{remark}
 To investigate the minimum generating set of $\mathcal{IC}_{n}$, we are going to  introduce the following definition.
 \begin{definition} A quasi-idempotent element  $\varepsilon$ of shift 1 as expressed in  \eqref{r} is called an \emph{essential element} if $y_{i}\varepsilon=y_{i}-1$, and so essential elements of height $p$ are of the form: \begin{equation}\label{r2}\varepsilon=\begin{pmatrix}y_1&\cdots&y_{i-1}&y_{i}&y_{i+1}&\cdots&y_p\\y_1&\cdots&y_{i-1}&y_{i}-1&y_{i+1}&\cdots&y_p\end{pmatrix} .\end{equation}
 \end{definition}

 \begin{remark} For $1\le n\le 2$, the quasi-idempotent elements in $\mathcal{IC}_{n}$ are all essential elements. However, for $n\geq 3$, not every quasi-idempotent element of shift 1 is an essential element. Moreover, every quasi-idempotent of height $n-1$ is necessarily essential.
 \end{remark}

 We now give and prove an enhanced version of [\cite{gmv}, Lemma 14.4.2] in the theorem below.

 \begin{theorem}\label{hq} Let $\varepsilon\in \mathcal{IC}_{n}$ be a quasi-idempotent element of shift 1 as expressed in \eqref{r}. Then $\varepsilon$  can be expressed as a product of essential elements.
 \end{theorem}
 \begin{proof} Let $\varepsilon\in \mathcal{IC}_{n}$ as expressed in \eqref{r} be a quasi-idempotent. Since $y_{i}>y_{i}\varepsilon$, we can let $s=y_{i}-y_{i}\varepsilon$, so clearly $y_{i}\varepsilon+j-1<y_{i}\varepsilon+j$ for all $1\leq j\le s$ and also $y_{i}=y_{i}\varepsilon+s$. Now for $1\leq j\le s$ define $\varepsilon_{{j}}$ as:
 \[\varepsilon_{{j}}=\begin{pmatrix}y_1&\cdots&y_{i-1}&y_{i}\varepsilon+j&y_{i+1}&\cdots&y_p\\y_1&\cdots&y_{i-1}&y_{i}\varepsilon+j-1&y_{i+1}&\cdots&y_p\end{pmatrix}.\]
\noindent Notice that  $y_{i-1}<y_{i}\varepsilon<y_{i}\varepsilon+1<\cdots< y_{i}\varepsilon+s=y_{i}<y_{i+1}$, and so for each $j$, we can easily see that $\varepsilon_{{j}}$ is an essential element. Moreover, it is not difficult to see that \[\varepsilon_{{s}} \varepsilon_{{s-1}}\cdots \varepsilon_{{1}}=\varepsilon.\]
\noindent Hence the result follows.
 \end{proof}

Notice that in the proof of  Lemma \ref{hq1} and Theorem \ref{hq}, $|\im \, \alpha|=h(\alpha)=h(\varepsilon)=h(\varepsilon_{{i}})=h(\varepsilon_{{j}})=p$ for all $1\le i\le p$ and all $1\le j\le s$. Therefore, we have the following result.

\begin{lemma}\label{hh}
Every element  in $S\in\{{RIC}_{n}(p),  \, K(n,p) \}$ of height $p$ can be expressed as a product of essential elements in $S$, each of height $p$.
\end{lemma}

We now have the following lemmas.

\begin{lemma}\label{es} For each $1\le p\le n-1$, there are $(n-1)\binom{n-2}{p-1}$ essential elements in ${RIC}_{n}(p)$.
\end{lemma}
\begin{proof} Notice that for $1\leq i\leq n-1$, the element $i+1$ can be paired with $i$ (which is equivalent to mapping $i+1$ to $i$)  in $n-1$ ways. Now, from the remaining $n-2$ elements (i.e., elements of $[n]\setminus\{i, i+1\}$), we can select  $p-1$ elements as our fixed points to form an essential element in $\binom{n-2}{p-1}$ ways. Hence we have altogether   $(n-1)\binom{n-2}{p-1}$ essential elements, and  the result now follows.
\end{proof}
It is worth noting that the number of essential elements in \( RIC_{n}(p) \) is represented by the triangle \(\underline{ A003506} \) in The On-Line Encyclopedia of Integer Sequences \cite{slone}. Now, let \( M(p) \) be the collection of all essential elements in \( RIC_{n}(p) \), and define \( G(p) = M(p) \cup E(RIC_{n}(p) \setminus \{0\}) \). Then, by Lemma \ref{es}, it follows that \( |M(p)| = (n-1) \binom{n-2}{p-1} \). We now have the following lemma.

\begin{lemma}\label{es3} For $1\le p\le n-1$,
$|G(p)|=(n-1)\binom{n-2}{p-1}+\binom{n}{p}$.
\end{lemma}
\begin{proof} The result follows directly from Lemma \ref{es} and Remark \ref{ls}(iv).
\end{proof}
The triangle of numbers from \( |G(p)| \) can be obtained by adding the sequences \( \underline{A003506} \) and \( \underline{A007318} \) in The On-Line Encyclopedia of Integer Sequences \cite{slone}, which is also equivalent to the sequence \( \underline{A103450} \), a figurate number triangle read by rows in The On-Line Encyclopedia of Integer Sequences \cite{slone}. The next result gives the total number of essential elements in the monoid \( \mathcal{IC}_{n} \). This number corresponds to the sequence \( \underline{A001787} \) in The On-Line Encyclopedia of Integer Sequences \cite{slone}.

\begin{theorem} Let $\mathcal{IC}_{n}$ be as defined in \eqref{qn111}. Then the number of essential elements in $\mathcal{IC}_{n}$ is $(n-1)\cdot2^{n-2}$.
\end{theorem}
\begin{proof} By summing over $p$ from $1$ to $n-1$ in Lemma \ref{es}, we see that:
\begin{align*} \sum\limits_{p=1}^{n-1}{(n-1)\binom{n-2}{p-1}}=& (n-1)\sum\limits_{p=1}^{n-1}\binom{n-2}{p-1}\\&= (n-1)\sum\limits_{k=0}^{n-2}{\binom{n-2}{k}} \, \, \left(\textnormal{by letting $k=p-1$, so that $p=k+1$, so  $0\leq k\leq n-2$.}\right)\\&=(n-1)\cdot 2^{n-2},
  \end{align*}
  \noindent as required.
\end{proof}

The next result shows that the set of nonzero idempotents together with the essential elements  in ${RIC}_{n}(p)$  is the minimum generating set of ${RICH}_{p}(n)$.
 \begin{lemma}\label{min4} Let $\alpha$, $\beta$ be elements in ${RIC}_{n}(p)\setminus \{0\}$ ($1\le p\le n$). Then $\alpha\beta\in G(p)$  if and only if $\alpha, \, \beta\in G(p)$  and $\alpha\beta=\alpha$  or $\alpha\beta=\beta$.
\end{lemma}

\begin{proof} Suppose
$\alpha\beta\in G(p)$. Thus either  $\alpha\beta\in E({RIC}_{n}(p)\setminus \{0\})$ or $\alpha\beta\in M(p)$. We consider the two cases separately.

\noindent \textbf{Case i.} Suppose $\alpha\beta\in E({RIC}_{n}(p)\setminus\{0\})$. Then
\[
p= f(\alpha\beta)\leq f(\alpha)\leq |\im \, \alpha|=p,\]
\[
p= f(\alpha\beta)\leq f(\beta)\leq |\im \, \beta|=p.\]

\noindent This ensures that \[F(\alpha)=F(\alpha\beta)=F(\beta), \] \noindent and so $\alpha, \, \beta\in E({RIC}_{n}(p)\setminus \{0\})$  and $\alpha\beta=\alpha$.

\noindent \textbf{Case ii.} Now suppose $\alpha\beta\in M(p)$. Thus $\alpha\beta$ is an essential element which can be expressed as in \eqref{r2}, that is: \[\alpha\beta=\begin{pmatrix}y_{1}& \cdots &y_{i-1}&y_{i}&y_{i-1}&\cdots&y_p\\y_{1}&\cdots&y_{i-1}&y_{i}-1&y_{i+1}&\cdots&y_p\end{pmatrix}.\] \noindent This means that $\dom \, \alpha=\dom \, \alpha\beta$, $\im \, \beta =\im \, \alpha\beta$ and  $\im \, \alpha=\dom \, \beta$. Thus,  \[\alpha=\begin{pmatrix}y_{1}& \cdots &y_{i-1}&y_{i}&y_{i+1}&\cdots&y_p\\y_{1}&\cdots&y_{i-1}&y_{i}\alpha&y_{i+1}&\cdots&y_p\end{pmatrix} \text{ and } \beta=\begin{pmatrix}y_{1}& \cdots &y_{i-1}&y_{i}\alpha&y_{i+1}&\cdots&y_p\\y_{1}&\cdots&y_{i-1}&y_{i}-1&y_{i+1}&\cdots&y_p\end{pmatrix}.\]  Notice that since $\alpha$ and $\beta$ are decreasing maps, we must have $y_{i}\alpha\leq y_{i}$ and $y_{i}-1\leq y_{i}\alpha$, and so  $y_{i}-1\leq y_{i}\alpha \leq y_{i}$. However, the fact that $\alpha\beta$ is an essential element, implies  either $y_{i}\alpha=y_{i}-1$ or $y_{i}\alpha=y_{i}$. We consider the two subcases separately.

\noindent \textbf{Subcase a.} If $y_{i}\alpha=y_{i}-1$, then $\alpha$ and $\beta$ are obviously essential and idempotent elements, respectively,  and so $\alpha\beta=\alpha$ and $\beta\in E({RIC}_{n}(p)\setminus \{0\})$.

\noindent \textbf{Subcase b.} If $y_{i}\alpha=y_{i}$, then $\alpha$ and $\beta$ are clearly idempotent and essential elements, respectively. Thus, it follows easily that  $\alpha\beta=\beta$ and $\alpha\in E({RIC}_{n}(p)\setminus \{0\})$. In either of the subcases, we see that $\alpha, \beta\in G(p)$ and either $\alpha\beta=\beta$ or $\alpha\beta=\alpha$. The converse is obvious.
\end{proof}

Thus, we have the following.

\begin{theorem}\label{pb} Let \( RIC_{n}(p) \) be as defined in \eqref{knn}. Then, the rank of \( RIC_{n}(p) \) is given by:
\[
\text{rank } RIC_{n}(p) = (n-1) \binom{n-2}{p-1} + \binom{n}{p}.
\]
\end{theorem}
\begin{proof} It follows from the fact that $G(p)$ is the minimum generating set, as stated by Lemma \ref{min4}, and its order is specified in Lemma \ref{es3}.
\end{proof}

The next lemma is crucial for determining the rank of the  injective partial Catalan monoid $\mathcal{IC}_{n}$.  Now,  for $1\leq p\leq n$, let \[J^{*}_{p}=\{\alpha\in \mathcal{IC}_{n}: |\im \, \alpha|=p \}.\]
\noindent Then we have:
\begin{lemma}\label{lm12} For $0\leq p\leq n-2$ and $n\geq 4$, $J^{*}_{p}\subset \langle J^{*}_{p+1}\rangle$. In other words, if $\alpha\in J^{*}_{p}$ then $\alpha\in \langle J^{*}_{p+1}\rangle$ for $1\leq p\leq n-2$.
\end{lemma}
\begin{proof} Using Theorem \ref{hq}, it suffices to prove that every element in \( G(p) \) can be expressed as a product of elements in \( G(p+1) \). That is to say, every idempotent of height \( p \) can be expressed as a product of idempotents of height \( p+1 \), and every essential element of height \( p \) can be expressed as a product of elements in \( G(p+1) \). Thus, we consider the elements of \( E(J^{*}_{p}) \) and \( M(p) \) separately.
\\

\noindent\textbf{i.} The  elements in $E(J^{*}_{p})$:\\

 Let $\epsilon\in E(J^{*}_{p})$ be expressed as: \[\epsilon= \begin{pmatrix}
y_1 &  \cdots  & y_p \\
y_1 &  \cdots  & y_p
\end{pmatrix},\] \noindent \noindent where $1\leq y_{1}<\cdots<y_{p}\leq n$. Since $p\leq n-2$, it follows that $(\dom \, \epsilon)^{\prime}$ contains at least two elements, say $c$ and $d$. Without loss of generality, suppose $c<d$. Let $A=\dom \, \epsilon \, \cup \, \{c\}$ and $B=\dom \, \epsilon \, \cup \, \{d\}$, and define $\epsilon_{1}$ and $\epsilon_{2}$ as follows:

For $x\in A$ and $y\in B$
\[x\epsilon_{1}=\left\{
                                                                                                                                \begin{array}{ll}
                                                                                                                                  x, & \hbox{if $x\neq c$;} \\
                                                                                                                                  c, & \hbox{if $x=c$}
                                                                                                                                \end{array}
                                                                                                                              \right.
  \quad \text{and} \quad
y\epsilon_{2}=\left\{
                                                                                                                                \begin{array}{ll}
                                                                                                                                  y, & \hbox{if $y\neq d$;} \\
                                                                                                                                  d, & \hbox{if $y=d.$}
                                                                                                                                \end{array}
                                                                                                                              \right.
\]

\noindent Clearly, $\epsilon_{1}$ and $\epsilon_{2}$ are idempotents in  $E(J^{*}_{p+1})$, and one can easily show that $\epsilon=\epsilon_{1}\epsilon_{2}$.

\noindent\textbf{ii.} The  elements in $M(p)$:\\

Let $\varepsilon$ be an essential element of height $p$ as expressed in \eqref{r2}, which has the form: \[\varepsilon=\begin{pmatrix}y_1&\cdots&y_{i-1}&y_{i}&y_{i+1}&\cdots&y_p\\y_1&\cdots&y_{i-1}&y_{i}-1&y_{i+1}&\cdots&y_p\end{pmatrix}.\]\noindent Now, since  $p\leq n-2$, it follows that $(\dom \, \varepsilon \, \cup \, \im \, \varepsilon)^{\prime}$ contains at least one element, say $d$. Notice that  $y_{i}-1\not\in \dom \, \varepsilon$.
Now let $A=\dom \, \varepsilon \, \cup \, \{d\}$ and $B=\dom \, \varepsilon \, \cup \, \{y_{i}-1\}$, and  define $\varepsilon^{\prime}$ and $\epsilon$ as follows:

For $x\in A$ and $y\in B$
\[x\varepsilon^{\prime}=\left\{
                                                                                                                                \begin{array}{ll}
                                                                                                                                  x\varepsilon, & \hbox{if $x\neq d $;} \\
                                                                                                                                  d & \hbox{if $x=d$}
                                                                                                                                \end{array}
                                                                                                                              \right.
 \quad \text{and} \quad
y\epsilon=y.\]
\noindent Notice that $\epsilon\in E(J^{*}_{p+1})\subset G(p+1)$, and it is not difficult to see that  $\varepsilon^{\prime}$ is an essential element in $M(p+1)\subset G(p+1)$. One can now easily show that $\varepsilon=\varepsilon^{\prime}\epsilon$. The proof of the lemma is now complete.
\end{proof}

It is important to note that that the rank of the two sided ideal  $K(n,p)$ has been found in [\cite{lu},  Theorem 2.1].
\begin{theorem}[\cite{lu},  Theorem 2.1]\label{knp} Let $K(n,p)$ be as defined in \eqref{kn}, where $1\leq p\leq n-1$. Then,  the rank of $K(n,p)$ is \[\text{rank } K(n,p)=(n-1)\binom{n-2}{p-1}+\binom{n}{p}.\]
\end{theorem}

Consequently, we obtain the following result.

\begin{theorem} Let $K(n,p)$ be as defined in \eqref{kn} and let \( RIC_{n}(p) \) be as defined in \eqref{knn}. Then the \(\text{rank } K(n,p) = \text{rank }RIC_{n}(p).\)
\end{theorem}
 We now deduce the following corollary to Theorem \ref{knp}.

\begin{corollary} The rank of \(\mathcal{IC}_{n}\) is:  \(\text{rank } \mathcal{IC}_{n}= 2n.\)
\end{corollary}
\begin{proof} Notice that by Lemma \ref{lm12}, $\langle J^{*}_{n-1} \rangle=\mathcal{IC}_{n}\setminus J^{*}_{n}=K(n,n-1)$. Also, observe that $J^{*}_{n}$ contains only the identity element $\text{id}_{[n]}$. Thus,   $\text{rank } \mathcal{IC}_{n}= \text{rank }K(n,n-1)+1$. The result now follows from Theorem \ref{knp}.
\end{proof}

\subsection{The rank of $\mathcal{Q}^{\prime}_{n}$}
In this subsection, we will generally consider an element \( \alpha \in \mathcal{Q}^{\prime}_{n} \) as expressed in \eqref{1}, with the property that \( 2 \leq x_{1} \), since the domains of elements in \( \mathcal{Q}^{\prime}_{n} \) do not contain 1. To investigate the rank of \( \mathcal{Q}^{\prime}_{n} \), we need the following definition and remark from \cite{zuf2}.

\begin{definition} An injective map  $\alpha$ in $ J^{*}_{p}=\{\alpha\in \mathcal{Q}^{\prime}_{n}: \, h(\alpha)=p\}$  is said  to be a \emph{requisite element} if it is of the form: \[\alpha_{i}=\begin{pmatrix}2&\cdots&i& a_{i}&\cdots& a_{p}\\1&\cdots&i-1&a_i&\cdots&a_{p}\end{pmatrix},\]
  \noindent where $1<i<a_{i}<a_{i+1}<\cdots< a_{p}\leq n$.
\end{definition}
\begin{remark} \label{requ}If $\alpha_{i}$ is a requisite element, then observe that for each $1\leq i\le p$, \[\dom \, \alpha_{i}=\{2,\ldots,i, a_{i}, \dots, a_{p}\}=\{2, \ldots,i\}\cup F(\alpha_{i}).\]
\noindent Moreover, $\alpha_{i}$ is unique in  $L^{*}_{\alpha_{i}}$, in the sense that no two $\alpha_{i}$'s belong to the same $\mathcal{L}^{*}-$class. However, an $\mathcal{R}^{*}-$class may contain more than one requisite element. In fact, in $J^{*}_{n-1}$ all the $n-1$ requisite elements belong to a single $\mathcal{R}^{*}-$class.
\end{remark}

 We immediately have the following lemma.
\begin{lemma}\label{ms} If $\alpha_{i}$ is the unique requisite element in $L^{*}_{\alpha}$  in $J_{p}^{*}$ \emph{(}$1\le p\le n-1$\emph{)}, then there exists $\beta\in R^{*}_{\alpha}$ such that $\alpha=\beta\alpha_{i}$.
\end{lemma}
 \begin{proof} Let \( \alpha_{i} \) be the unique requisite element in \( L^{*}_{\alpha} \) in \( J^{*}_{p} \), where \( \alpha \) is expressed as in \eqref{1}, with \( 2 \leq x_{i} \). Now, since \( \alpha_{i} \in L^{*}_{\alpha} \), it follows that \( \im \, \alpha = \im \, \alpha_{i} \), and thus \( a_{j} = j \) for all \( 1 \leq j \leq i - 1 \). Therefore, \( \alpha \) and \( \alpha_{i} \) are given by:
 \[\alpha=\begin{pmatrix}x_{1}&\cdots&x_{i-1}& x_{i}&\cdots& x_{p}\\1&\cdots&i-1&a_i&\cdots&a_{p}\end{pmatrix} \text{ and }\alpha_{i}=\begin{pmatrix}2&\cdots&i& a_{i}&\cdots& a_{p}\\1&\cdots&i-1&a_i&\cdots&a_{p}\end{pmatrix},\]\noindent respectively.
By the isotone property, it is clear that $x_{j}<x_{j+1}$ for all $1\le j\le i-1$, and since $2\leq x_{1}$, it follows that $j+1\leq x_{j}$. Thus,  the map $\beta$ defined as:
\[\beta=\begin{pmatrix}x_{1}&\cdots&x_{i-1}& x_{i}&\cdots& x_{p}\\2&\cdots&i&a_i&\cdots&a_{p}\end{pmatrix}\in \mathcal{Q}^{\prime}_{n}.\]
\noindent  Furthermore, notice that $\dom \, \alpha=\dom \, \beta$, and so by Remark \ref{a5}, we see that $\beta\in R^{*}_{\alpha}$. Clearly, $\beta\alpha_{i}=\alpha$. The proof is now complete.
\end{proof}


Next, we prove the following theorem.
\begin{theorem}\label{hqa}Let \( \alpha \in \mathcal{Q}^{\prime}_{n} \) be expressed as in \eqref{1}, where \( 2 \leq x_{1} \). Then \begin{itemize}
                                                                            \item[(i)] if $a_{1}\neq 1$, then $\alpha$ is generated by essential elements;
                                                                            \item[(ii)] if $a_{1}= 1$, then $\alpha$ is a product of essential elements and the unique requisite element in $L^{*}_{\alpha}$.
                                                                          \end{itemize}
\end{theorem}
\begin{proof} Let $\alpha\in \mathcal{Q}^{\prime}_{n}$ be as expressed in  \eqref{1}, where $2\leq x_{1}$.

\noindent\textbf{(i.)} Suppose $a_{1}\neq 1$ and let $U=\{\alpha\in\mathcal{Q}^{\prime}_{n}: \, 1\not\in \, \im \, \alpha \}$. Then it is not difficult to see that $U$ is a subsemigroup of $\mathcal{Q}^{\prime}_{n}$, which is isomorphic to $\mathcal{IC}_{n-1}$. Hence by Theorem \ref{hq}, $U$ is  generated by essential elements.

 \noindent\textbf{ (ii.)} Now suppose $a_{1}=1$. Thus, by Lemma \ref{ms}, $\alpha$ can be expressed as  \[\alpha=\beta\alpha_{i},\] \noindent for some  $\beta\in R^{*}_{\alpha}$, where $\alpha_{i}$ is the unique requisite  element in $L^{*}_{\alpha}$. To be precise,
 \[\beta=\begin{pmatrix}x_{1}&\cdots&x_{i-1}& x_{i}&\cdots& x_{p}\\2&\cdots&i&a_i&\cdots&a_{p}\end{pmatrix} \text{ and }\alpha_{i}=\begin{pmatrix}2&\cdots&i& a_{i}&\cdots& a_{p}\\1&\cdots&i-1&a_i&\cdots&a_{p}\end{pmatrix}.\]
\noindent Notice that  $1\not\in\im \, \beta$, and so $\beta\in U$. Thus, by (i) $\beta$ is generated by  essential elements. The result now follows.
\end{proof}

We now have the following corollary.
\begin{corollary} The semigroup $\mathcal{Q}^{\prime}_{n}$ is generated by its idempotents,  essential and  requisite elements.
\end{corollary}
Notice that in the proof of Theorem \ref{hqa}, $|\im \, \alpha|=h(\alpha)=h(\beta)=h(\alpha_{i})=p$ for all $1\le i\le p$. Based on this, we have the following result.

\begin{lemma}\label{hh}
Every element  in $\mathcal{Q}^{\prime}_{n}$  of height $p$ can be expressed as a product of idempotents, essential and  requisite elements in $\mathcal{Q}^{\prime}_{n}$, each of height $p$.
\end{lemma}

 Now let $R(p)=\{\alpha\in {RQ}^{\prime}_{n}(p): \alpha \text{ is a requisite element}\}$ and $M(p)=\{\alpha\in {RQ}^{\prime}_{n}(p): \alpha \text{ is an essential element} \}$. Then we have the following lemma.

  \begin{lemma}\label{nid}  For $1\le p\leq n-1$, we have \begin{itemize} \item[(i)] $|R(p)|=\binom{n-1}{p-1}$; \item[(ii)] $|E({RQ}^{\prime}_{p}(n)\setminus\{0\})|={\binom{n-1}{p}}$;
  \item[(iii)]$|M(p)|=(n-2)\binom{n-3}{p-1}$  \emph{(}$1\le p\leq n-2$\emph{)}.
  \end{itemize}
\end{lemma}
\begin{proof}
\noindent \textbf{(i.)} Given that $1$ must be included in every image set that contains the $p$ images, we can then choose the remaining $p-1$ images from the available $n-1$ elements in $\binom{n-1}{p-1}$ distinct ways. However, notice that the domain of  any requisite element is uniquely chosen once the image set is fixed.

\noindent\textbf{(ii.)} Notice that since $1$ is not in the domain of every element in ${RQ}^{\prime}_{p}(n)\setminus\{0\}$, then to count the idempotent of height $p$, is equivalent to counting all possible subsets of $[n]\setminus\{1\}$ with $p$ elements. This is equivalent  to selecting $p$ elements from $n-1$ elements.

\noindent\textbf{(iii.)} Let $\epsilon$ be an essential  as expressed in \eqref{r2}. Notice that for \( 2 \leq i \leq n - 1 \), the element \( i + 1 \) can be paired with \( i \) in \( n - 2 \) ways. Now, from the  \( n - 3 \) remaining elements (i.e., elements of \( [n] \setminus \{1, i, i + 1\} \)), we can select the \( p - 1 \) fixed points  in \( \binom{n - 3}{p - 1} \) ways. Hence, we have altogether \( (n - 2) \binom{n - 3}{p - 1} \) essential elements, and  the result now follows.
\end{proof}

Now, let $G(p)=(E({RQ}^{\prime}_{p}(n)\setminus\{0\}))\cup M(p)\cup R(p)$. Then as a direct consequence of the above lemma we have the following lemma.
 \begin{lemma}\label{og} If $G(p)=E({RQ}^{\prime}_{n}(p)\setminus\{0\})\cup M(p)\cup R(p)$, then $$|G(p)|=\binom{n-1}{p-1}+{\binom{n-1}{p}}+(n-2)\binom{n-3}{p-1}=\binom{n}{p}+(n-2)\binom{n-3}{p-1}.$$
 \end{lemma}

 The next result shows that the subset $G(p)$ of ${RQ}^{\prime}_{n}(p)$  is the minimum generating set of ${RQ}^{\prime}_{n}(p)$.
 \begin{lemma}\label{minnn} Let $\alpha$, $\beta$ be elements in ${RQ}^{\prime}_{n}(p)\setminus \{0\}$. Then $\alpha\beta\in G(p)$  if and only if $\alpha, \, \beta\in G(p)$  and $\alpha\beta=\alpha$ or $\alpha\beta=\beta$.
\end{lemma}

\begin{proof} Suppose $\alpha\beta\in G(p)$. Thus, either \textbf{(i.)} $\alpha\beta\in E({RQ}^{\prime}_{n}(p)\setminus\{0\})$ or \textbf{(ii.)} $\alpha\beta\in R(p)$ or \textbf{(iii.)} $\alpha\beta\in M(p).$ We consider the three cases separately.

Notice that since the idempotent and essential elements in \( RQ^{\prime}_{n}(p) \setminus \{0\} \) are also idempotent and essential elements in \( RIC^{\prime}_{p}(n-1) \setminus \{0\} \) (by isomorphism), then \textbf{(i)} and \textbf{(ii)} follow from Lemma \ref{min4}.

\noindent \textbf{ (iii.)} Now suppose $\alpha\beta\in M(p)$. Thus $\alpha\beta$ is a requisite element which has the form \[\alpha\beta=\begin{pmatrix}2& \cdots &i&a_{i}&\cdots&a_p\\1&\cdots&i-1&a_{i}&\cdots&a_p\end{pmatrix},\] \noindent where $1<i<a_{i}<a_{i+1}<\cdots< a_{p}\leq n$. This means that $\dom \, \alpha=\dom \, \alpha\beta$, $\im \, \beta =\im \, \alpha\beta$ and  $\im \, \alpha=\dom \, \beta$. Thus, \[\alpha=\begin{pmatrix}2& \cdots &i&a_{i}&\cdots&a_p\\2\alpha&\cdots&i\alpha&a_{i}&\cdots&a_p\end{pmatrix} \text{ and } \beta=\begin{pmatrix}2\alpha& \cdots &i\alpha&a_{i}&\cdots&a_p\\1&\cdots&i-1&a_{i}&\cdots&a_p\end{pmatrix}.\]
\noindent The claim here is that $\alpha$ must be an idempotent. Notice that either $2\alpha=1$ or $2\alpha=2$. In the former, we see that $1=2\alpha\in \dom \, \beta$ which is a contradiction. Therefore $2\alpha=2$ which implies $j\alpha=j$ for all $2\leq j\leq i$. Thus $\alpha$ is an idempotent, that is, \[\alpha=\begin{pmatrix}2& \cdots &i&a_{i}&\cdots&a_p\\2&\cdots&i&a_{i}&\cdots&a_p\end{pmatrix} \text{ and }\beta=\begin{pmatrix}2& \cdots &i&a_{i}&\cdots&a_p\\1&\cdots&i-1&a_{i}&\cdots&a_p\end{pmatrix}.\] \noindent  Therefore, $\beta\in M(p)\subset G(p)$ and $\alpha\in E({RQ}^{\prime}_{p}(n)\setminus\{0\})\subset G(p)$, and also $\alpha\beta=\beta.$
The converse is obvious.
\end{proof}

At this point, we present a crucial result in this section.

\begin{theorem}\label{pbb} Let ${RQ}^{\prime}_{n}(p)$ be as defined in \eqref{mnnn}. Then, the rank of $ {RQ}^{\prime}_{n}(p)$ is:  \[\text{rank } {RQ}^{\prime}_{n}(p)= \binom{n}{p}+(n-2)\binom{n-3}{p-1}.\]
\end{theorem}
\begin{proof} The proof follows from Lemmas \ref{minnn} \& \ref{og}.
\end{proof}

The next lemma is crucial in determining  the ranks of the  semigroup  $\mathcal{Q}^{\prime}_{n}$ and its two sided ideal $M(n,p)$.  Now,  for $1\leq p\leq n-1$, let \[J^{*}_{p}=\{\alpha\in \mathcal{Q}^{\prime}_{n}: |\im \, \alpha|=p \}.\] Moreover, for $1\leq p\leq n-1$ let $M(p)$ be the collection of all essential elements in $J_{p}^{*}$, and let $R(p)$ be the collection of all  requisite elements in $J_{p}^{*}$, and let \(G(p)=M(p)\cup E(J^{*}_{p}).\) Then we have the following lemmas.
 \begin{lemma} For $1\le p\le n-1$, $|G(p)|=\binom{n}{p}+(n-2)\binom{n-3}{p-1}$.
\end{lemma}
\begin{proof} The result follows from Lemma \ref{og}.
\end{proof}
\begin{lemma}\label{lm1} For $0\leq p\leq n-3$, $J^{*}_{p}\subset \langle J^{*}_{p+1}\rangle$. In other words, if $\alpha\in J^{*}_{p}$ then $\alpha\in \langle J^{*}_{p+1}\rangle$ for $1\leq p\leq n-3$.
\end{lemma}
\begin{proof} Using Theorem \ref{hqa}, it suffices to prove that every element in \( G(p) \) can be expressed as a product of elements in \( G(p+1) \). That is to say, \textbf{(i)} every idempotent of height \( p \) can be expressed as a product of idempotents, each of height \( p + 1 \); \textbf{(ii)} every essential element of height \( p \) can be expressed as a product of essential elements, each of height \( p + 1 \); and \textbf{(iii)} every requisite element can be expressed as a product of requisite elements, each of height \( p + 1 \). However, notice that if \( \epsilon \in E(J^{*}_{p}) \) and \( \alpha \in M(p) \), then \( \epsilon \) and \( \alpha \) will be an idempotent and an essential element in \( \mathcal{IC}_{n-1} \) (by isomorphism). Therefore, \textbf{(i)} and \textbf{(ii)} follow from Lemma \ref{lm1}.\\

\noindent\textbf{(iii.)} The  elements in $M(p)$:  Let $\alpha_{i}$ be a requisite element of height $p$, which has the form: \[\alpha_{i}= \begin{pmatrix}2& \cdots &i&a_{i}&\cdots&a_p\\1&\cdots&i-1&a_{i}&\cdots&a_p\end{pmatrix},  \]\noindent where  $1<i<a_{i}<a_{i+1}<\cdots< a_{p}\leq n$.
Now, since  $p\leq n-3$, it implies that $(\dom \, \alpha_{i} \, \cup  \, \im \, \alpha_{i})^{c}$ contains at least two elements, specifically $c$ and $d$. We may suppose without loss of generality that, $d<c$.
Now let $A=\dom \, \alpha_{i} \, \cup \, \{d\}$ and $B=\dom \, \alpha_{i} \, \cup \, \{c\}$. Now  define $\beta$ and $\gamma$ as follows:

For $x\in A$ and $y\in B$
\[x\beta=\left\{
                                                                                                                                \begin{array}{ll}
                                                                                                                                  x, & \hbox{if $x\neq d $;} \\
                                                                                                                                  d, & \hbox{if $x=d$}
                                                                                                                                \end{array}
                                                                                                                              \right.
\quad \text{and} \quad
y\gamma=\left\{
                                                                                                                                \begin{array}{ll}
                                                                                                                                  y\alpha_{i}, & \hbox{if $y\neq c$;} \\
                                                                                                                                  c, & \hbox{if $y=c$.}
                                                                                                                                \end{array}
                                                                                                                              \right.
\]
\noindent Notice that $\beta\in E(J^{*}_{p+1})\subset G(p+1)$, and it is not difficult to see that  $\gamma$ is a requisite element in $M(p+1)\subset G(p+1)$. One can now easily show that $\alpha_{i}=\beta\gamma$. The proof of the lemma is now complete.
\end{proof}

Therefore, we obtain  have the following result.

\begin{theorem}\label{knp2} Let $M(n,p)$ be as defined in \eqref{mn}. Then  for $1\leq p\leq n-2$, we have \[\text{rank } M(n,p) =\binom{n}{p}+(n-2)\binom{n-3}{p-1}.\]
\end{theorem}
\begin{proof}Notice that by Lemma \ref{lm1},   $\langle J^{*}_{p} \rangle= M(n,p)$ for all $1\le p\leq n-2$.  Notice also that, $\langle E({RQ}^{\prime}_{n}(p)\setminus\{0\})\cup R(p)\cup M(p)\rangle= J^{*}_{p}$. The result now follows from Theorem \ref{pbb}.
\end{proof}

It is important to note that Lemma \ref{lm1} does not cover the case \( p = n - 2 \), meaning that the assertion \( J^{*}_{n-2} \subset \langle J^{*}_{n-1} \rangle \) is not true. This can be demonstrated by first noting that there are some essential elements in \( J^{*}_{n-2} \) which cannot be  generated by  elements of \( J^{*}_{n-1} \). We will demonstrate this by first writing the form of all the essential elements in \( J^{*}_{n-2} \) and the requisite elements in \( J^{*}_{n-1} \). Now, for \( 2 \leq i \leq n - 1 \), the essential elements in \( J^{*}_{n-2} \) are of the form:
\begin{equation}\label{11}
\alpha_{i,i+1} = \begin{pmatrix}
2 & \cdots & i - 1 & i + 1 & i + 2 & \cdots & n \\
2 & \cdots & i - 1 & i & i + 2 & \cdots & n
\end{pmatrix}.
\end{equation}
\noindent Notice that \( i \notin \dom\,\alpha_{i,i+1} \) and \( i + 1 \notin \im\,\alpha_{i,i+1} \).

For $2\leq i\leq n$, the requisite elements in $J^{*}_{n-1}$ are of the form  \begin{equation}\label{kk-1}\alpha_{i}=\begin{pmatrix}2&\cdots&i-1& {i}&i+1&\cdots& n\\1&\cdots&i-2&i-1&i+1&\cdots&n\end{pmatrix}.\end{equation}
\noindent Notice also that $i\not\in \im \, \alpha_{i}$. The unique idempotent  \begin{equation}\label{13}\epsilon=\left(\begin{array}{cccccccc}
                                                                            {2}&  \cdots &  n \\
                                                                            {2} & \cdots &   n
                                                                          \end{array}
\right)\in J^{*}_{n-1} ,\end{equation}
\noindent  is the unique \emph{left identity} in \( \mathcal{Q}^{\prime}_{n} \), but not an identity. Perhaps we may call the semigroup \( \mathcal{Q}^{\prime}_{n} \) a \emph{right ample left monoid}. It is also clear that \( J^{*}_{n-1} \) has only one \( \mathcal{R}^{*} \)-class and \( n \) \( \mathcal{L}^{*} \)-classes. Obviously, there are \( n - 1 \) \( \mathcal{L}^{*} \)-classes, each containing a unique requisite element.

Now, it is clear that \( 2 \notin \dom\, {\alpha}_{ i}\epsilon\cup \dom\, \alpha_{i}\alpha_{j} \)  for all $i$ and  $j$. Thus, all essential elements  in \( J^{*}_{n-2} \) each of whose domain contains the element 2 cannot be generated by  elements in \( J^{*}_{n-1} \). However,  one can easily see that $\alpha_{3}^2={\alpha}_{ 2,3}$. That is to say the essential element ${\alpha}_{ 2,3}$ in $J^{*}_{n-2}$ can be generated by the essential element $\alpha_{3}$ in \( J^{*}_{n-1} \). Hence the set of essential elements   $M(n-2)\setminus\{{\alpha}_{ 2,3}\}$ in $J^{*}_{n-2}$ is enough to be part of any minimal generating set, where ${\alpha}_{ 2,3}=\left(\begin{array}{cccc}
                                                                            {3}&4&  \cdots& n \\
                                                                            {2} & 4 &\cdots& n
                                                                          \end{array}
\right)$.

It is clear that $G(n-1)=R(n-1)\cup\{\epsilon\}=J^{*}_{n-1}$, and so $|G(n-1)|=n$.  We now have the following:
\begin{lemma}\label{minnnn} Let $\alpha$, $\beta$ be elements in $\langle J^{*}_{n-1}\rangle$ . Then $\alpha\beta\in G(n-1)$  if and only if $\alpha, \, \beta\in G(n-1)$  and $\alpha\beta=\alpha$ or $\alpha\beta=\beta$.
\end{lemma}
\begin{proof} The proof is similar to that of Lemma \ref{minnn}.
\end{proof}
\begin{lemma}\label{jnn} The rank$\langle J^{*}_{n-1} \rangle=n$.
\end{lemma}
\begin{proof} The result follows from Lemma \ref{minnnn} and the fact that $|G(n-1)|=|J^{*}_{n-1}|$.
\end{proof}

 The requisite elements in $J^*_{n-2}$ are generally of the form: \begin{equation}\label{iki}{\alpha}_{i, k}=\left(\begin{array}{ccccccccc}
                                                                            {2}& \cdots &i &i+1& \cdots& k-1 &   k+1 & \cdots& n \\
                                                                            {1} & \cdots & i-1 & i+1 &\cdots& k-1  &k+1 &\cdots & n
                                                                          \end{array}
\right)\quad (2\leq i< k\leq n).\end{equation}\noindent Notice that $1, \, k\not\in \dom \, {\alpha}_{ i, k}$ and $i,k\not\in \im \, {\alpha}_{[k, i]}$. Also, the element

\begin{equation}\label{i1i}{\epsilon}_{ 1,k}=\left(\begin{array}{ccccccc}
                                                                            {2}&\cdots &k-1 &k+1& \cdots& n \\
                                                                            {2} & \cdots & k-1 & k+1 &\cdots& n
                                                                          \end{array}
\right)\end{equation} \noindent is a partial identity in $G(n-2)$ for all $2\leq k\leq n$. Notice also  that $1, \, k\not\in \dom \, {\epsilon}_{  1, k}$ and $1,k\not\in \im \, {\epsilon}_{  1, k}$.

The subsequent lemma demonstrates that the requisite elements in $G(n-2)$ are not required to be included in any minimal  generating set of $\langle J^{*}_{n-2} \cup J^{*}_{n-1} \rangle$. However, the requisite elements in $G(n-1)$ suffice to be included in any minimal generating set of $\langle J^{*}_{n-2} \cup J^{*}_{n-1} \rangle$.
\begin{lemma}\label{ss1} $\langle G(n-2) \rangle \subseteq \langle \left(G(n-2)\setminus (R(n-2)\cup\{\alpha_{2,3}\})\right)\cup G(n-1) \rangle=\langle J^{*}_{n-2}\cup J^{*}_{n-1}  \rangle=\mathcal{Q}^{\prime}_{n}$.
\end{lemma}

\begin{proof}To show the inclusion, it is enough to show that $(R(n-2)\cup \{\alpha_{2,3}\})\subset \langle G(n-1)\cup E(J^{*}_{n-2}) \rangle$.
 Now  let   ${\alpha}_{  i, k} \in R(n-2)$ be  as expressed in \eqref{iki}. Next, take  ${\epsilon}_{1,k}\in E(J^{*}_{n-2})$ and   $\alpha_{i}\in G(n-1)$  as expressed in \eqref{i1i} \& \eqref{kk-1}, respectively; and observe that for any $i<k$, we have: \begin{align*}{\epsilon}_{ 1,k}\alpha_{i}&=\left(\begin{array}{ccccccccccc}
                                                                            {2}&\cdots &i&i+1&\cdots&k-1 &k+1& \cdots& n \\
                                                                            {2} & \cdots &i&i+1&\cdots& k-1 & k+1 &\cdots& n
                                                                          \end{array}
\right)\left(\begin{array}{cccccccc}
                                                                            {2}& 3& \cdots  &i& i+1&\cdots&  n \\
                                                                            {1} & 2&\cdots & i-1 & i+1&\cdots&   n
                                                                          \end{array}
\right)\\&=\left(\begin{array}{ccccccccccc}
                                                                            {2}&\cdots& i&i+1&\cdots&k-1 &k+1& \cdots& n \\
                                                                            {1} & \cdots &i-1&i+1&\cdots& k-1 & k+1 &\cdots& n
                                                                          \end{array}
\right)={\alpha}_{ i,k}.
\end{align*}
\noindent Moreover, \begin{align*}\label{2,3}{\alpha}_{ 3}\epsilon&=\left(\begin{array}{ccccccccccc}
                                                                            {2}&3 &4& \cdots& n \\
                                                                            {1} & 2 &4&\cdots& n
                                                                          \end{array}
\right)\left(\begin{array}{cccccccc}
                                                                            {2}& \cdots&  n \\
                                                                            {2} & \cdots&   n
                                                                          \end{array}
\right)=\left(\begin{array}{ccccccccccc}
                                                                            {3}&4&  \cdots& n \\
                                                                            {2} & 4 &\cdots& n
                                                                          \end{array}
\right)={\alpha}_{ 2,3},
\end{align*}\noindent where ${\alpha}_{ 3}, \, \epsilon\in G(n-1)$.

The equality is obvious. Hence, the result now follows.
\end{proof}
Now let $W=\left(G(n-2)\setminus (R(n-2)\cup\{\alpha_{2,3}\})\right)\cup G(n-1)$, then we have the following lemma.

\begin{lemma}\label{minnnnn} Let $\alpha, \beta\in J^{*}_{n-2}\cup J^{*}_{n-1}$. Then $\alpha\beta\in W$ if and only if $\alpha, \beta\in W$ and $\alpha\beta=\alpha$ or $\alpha\beta=\beta$.
\end{lemma}
\begin{proof} The proof follows from Lemmas \ref{minnn},  \ref{minnnn} \& \ref{ss1}.
\end{proof}
Finally, we conclude this section with the following result:

\begin{theorem} Let $\mathcal{Q}_{n}^{\prime}$ be as defined in \eqref{qn2}. Then for $n>1$, the rank $\mathcal{Q}_{n}^{\prime}=n^{2}-3n+4$.
\end{theorem}
\begin{proof} Clearly by Lemma \ref{minnnnn}, $\left(G(n-2)\setminus (R(n-2)\cup \{\alpha_{2,3}\})\right)\cup G(n-1)$ is the minimum generating set of $\langle J^{*}_{n-2}\cup J^{*}_{n-1} \rangle=\mathcal{Q}_{n}^{\prime}$ , and so by Lemmas \ref{nid} \&  \ref{jnn} we see that
\begin{align*} \text{ rank }\mathcal{Q}_{n}^{\prime}&= |G(n-2)|- (|R(n-2)|-1)+|G(n-1)|\\&=\binom{n-1}{n-3}+\binom{n-1}{n-2}+(n-2)\binom{n-3}{n-4}-\binom{n-1}{n-3}+n-1 \\&=\binom{n-1}{n-2}+(n-2)\binom{n-3}{n-4}+n-1\\&=(n-1)+(n-2)(n-3)+n-1=n^{2}-3n+4,
\end{align*}
as postulated.
\end{proof}

\section{The maximal subsemigroups of $\mathcal{IC}_{n}$}
A subsemigroup $ M \subseteq S $ is called \emph{maximal} provided that $ M \neq S $ and for any subsemigroup $ X \subseteq S $, the inclusion $ M \subseteq X $ implies $ M = X $ or $ X = S $. That is, a proper non-empty subsemigroup $M$ of a semigroup $S$, is maximal if $M\subseteq T\subseteq S$ for some subsemigroup $T$ of $S$, we have $M=T$ or $T=S$ \cite{gmv, gra}. In other words, a proper subsemigroup of a semigroup $S$ is considered maximal if it is not contained in any other proper subsemigroup of $S$. We would like to draw the reader's attention to [\cite{east}, Table 1], where a list of various transformation semigroups and their  number of  maximal subsemigroups  is provided. An element $a\in S$ is said to be \emph{indecomposable} if there are no elements $b, c\in S\setminus\{a\}$ with $b\neq c$ such that $a=bc$.

 We initiate our findings by first  recalling from Lemma \ref{min4} that the set $G(n-1)\cup\{\text{id}_{[n]}\}$ is the minimum generating set of $\mathcal{IC}_{n}$, as such an element $\alpha\in \mathcal{IC}_{n} $  is indecomposable if and only if  $\alpha\in G(n-1)\cup\{\text{id}_{[n]}\}$. Recall also that on the semigroup $\mathcal{IC}_{n}$, the set $G(n-1)$ consists of the essential and the idempotent elements of $J^{*}_{n-1}$. The essential elements and idempotents are of the form \begin{equation}\label{ele1} \varepsilon_{i,i+1}=\left(\begin{array}{cccccccc}
                                                                            {1}&  \cdots  &i-1& i+1&i+2&\cdots&  n \\
                                                                            {1} &\cdots & i-1 & i&i+2&\cdots&   n
                                                                          \end{array}
\right) \quad(1\le i\le n-1), \quad \text{and } \varepsilon_{i}=\text{id}_{[n]\setminus\{i\}} \quad (1\le i\le n),\end{equation}
\noindent respectively.
Thus, we have the following result.
\begin{theorem}\label{maxim} Let $\varepsilon_{i,i+1}$ and $\varepsilon_{i}$ be as defined in  \eqref{ele1}. A subsemigroup $M$ of $\mathcal{{IC}}_n$ is maximal if and only if $M$ belongs to one of the following three types:

\begin{itemize}
\item[\textnormal{(}i\textnormal{)}]  $M_{\textit{id}_{[n]}} := \mathcal{{IC}}_n \symbol{92} {\{\textit{id}_{[n]}\}};$
\item[\textnormal{(}ii\textnormal{)}] $M_{\varepsilon_{i,i+1}} := \mathcal{{IC}}_n \symbol{92} \{\varepsilon_{i,i+1}\}, $ for each $1\le i\le n-1$;
\item[\textnormal{(}iii\textnormal{)}] $M_{\varepsilon_{i}} := \mathcal{{IC}}_n \symbol{92} \{\varepsilon_{i}\},$ for each  $1\le i\le n$.
\end{itemize}
\end{theorem}
\begin{proof}
The proof follows from the fact that \( G(n-1) \) is the minimum generating set consisting of all the indecomposable elements of \( \mathcal{IC}_n \). Specifically, \( G(n-1) \) contains the idempotents \( \varepsilon_{i} \), the essential elements \( \varepsilon_{i,i+1} \), and the identity element \( \text{id}_{[n]} \).
\end{proof}
Therefore, we have the following.
\begin{corollary} The injective partial Catalan monoid $\mathcal{{IC}}_n$ contains exactly $2n$ maximal subsemigroups.
\end{corollary}
\begin{proof} The result follows from counting all the maximal subsemigroups as stated in Theorem \ref{maxim}.
\end{proof}

In conclusion, we wrap up the paper with the following results without proof.

\begin{theorem}\label{maxim} Let $\alpha_{i,i+1}$, $\alpha_{i}$ and $\epsilon$ be as defined in  \eqref{11}, \eqref{kk-1} and \eqref{13}, respectively. A subsemigroup $M$ of $\mathcal{{Q}}_n^{\prime}$ is maximal if and only if $M$ belongs to one of the following three types:

\begin{itemize}
\item[\textnormal{(}i\textnormal{)}]  $M_{\epsilon} := \mathcal{{Q}}_{n}^{\prime} \symbol{92} {\{\epsilon\}};$
\item[\textnormal{(}ii\textnormal{)}] $M_{\alpha_{i,i+1}} := \mathcal{{Q}}_{n}^{\prime} \symbol{92} \{\alpha_{i,i+1}\}, $ for each $i\in\{3,4,\ldots,    n-1\}$;
\item[\textnormal{(}iii\textnormal{)}] $M_{\alpha_{i}} := \mathcal{{Q}}_{n}^{\prime} \symbol{92} \{\alpha_{i}\},$ for each  $2\le i\le n$.
\end{itemize}
\end{theorem}
\begin{corollary} The semigroup $\mathcal{{Q}}_n^{\prime}$ contains exactly $n^{2}-3n+4$ maximal subsemigroups.
\end{corollary}

\noindent{\bf Acknowledgements, Funding and/or Conflicts of interests/Competing interests.} The third named author would like to thank Bayero University and
TETFund (TETF/ES/UNI/KANO/TSAS/2022) for financial support. He would also like to thank Sultan Qaboos University, Oman,  for hospitality during a 1-year postdoctoral research visit to the institution.

\end{document}